\documentclass[british]
{elsarticle}
\usepackage{amsmath, amssymb} 
\usepackage{
eucal,url,
enumerate,amscd,
}

\begin{document}


\title
{On the Geometry of Connections with Totally Skew-Symmetric
Torsion on Manifolds with Additional Tensor Structures and
Indefinite Metrics}

\author{Mancho Manev}
\ead{mmanev@uni-plovdiv.bg}\ead[url]{http://fmi.uni-plovdiv.bg/manev/}
\author{Dimitar Mekerov}
\ead{mircho@uni-plovdiv.bg}
\author{Kostadin Gribachev}
\ead{costas@uni-plovdiv.bg}
\address{Paisii Hilendarski University of Plovdiv \\
Faculty of Mathematics and Informatics\\
236 Bulgaria Blvd., 4003 Plovdiv, Bulgaria}




\newcommand{\ie}{i.~e. }
\newcommand{\X}{\mathfrak{X}}
\newcommand{\W}{\mathcal{W}}
\newcommand{\K}{\mathcal{K}}
\newcommand{\TT}{\mathcal{T}}
\newcommand{\N}{\mathbb{N}}
\newcommand{\R}{\mathbb{R}}
\newcommand{\s}{\mathfrak{S}}
\newcommand{\n}{\nabla}
\newcommand{\D}{{\rm d}}
\newcommand{\tr}{\mathrm{tr}}
\newcommand{\al}{\alpha}
\newcommand{\bt}{\beta}
\newcommand{\gm}{\gamma}
\newcommand{\lm}{\lambda}
\newcommand{\ea}{\varepsilon_\alpha}
\newcommand{\eb}{\varepsilon_\beta}
\newcommand{\eg}{\varepsilon_\gamma}
\newcommand{\sa}{\sum_{\al=1}^3}
\newcommand{\sbt}{\sum_{\bt=1}^3}
\newcommand{\ee}{\end{equation}}
\newcommand{\be}[1]{\begin{equation}\label{#1}}
\def\bea{\begin{eqnarray}} \def\eea{\end{eqnarray}}
\newcommand{\norm}[1]{\left\Vert#1\right\Vert ^2}
\newcommand{\nJ}[1]{\norm{\nabla J_{#1}}}
\newcommand{\Id}{\textrm{Id}}
\newcommand{\nT}{\norm{T}}

\newcommand{\ff}{\varphi}
\newcommand{\F}{\mathcal{F}}
\newcommand{\G}{\mathcal{G}}
\newcommand{\sx}{\mathop{\mathfrak{S}}\limits_{x,y,z}}
\newcommand{\nf}{\norm{\n \ff}}
\newcommand{\M}{(M,\ff,\xi,\eta,g)}
\newcommand{\Lf}{(G,\ff,\xi,\eta,g)}
\newcommand{\dd}{\textrm{d}}

\newcommand{\thmref}[1]{The\-o\-rem~\ref{#1}}
\newcommand{\propref}[1]{Pro\-po\-si\-ti\-on~\ref{#1}}
\newcommand{\secref}[1]{\S\ref{#1}}
\newcommand{\lemref}[1]{Lem\-ma~\ref{#1}}
\newcommand{\dfnref}[1]{De\-fi\-ni\-ti\-on~\ref{#1}}
\newcommand{\coref}[1]{Corollary~\ref{#1}}


\newtheorem{thm}{Theorem}
\newtheorem{lem}[thm]{Lemma}
\newtheorem{prop}[thm]{Proposition}
\newtheorem{cor}[thm]{Corollary}
\newdefinition{rmk}{Remark}
\newdefinition{dfn}{Definition}
\newdefinition{ack}{Acknowledgements}
\newproof{pf}{Proof}
\newproof{pot}{Proof of Theorem \ref{thm-geom}}




\hyphenation{Her-mi-ti-an ma-ni-fold ah-ler-ian}



\begin{abstract}
This paper is a survey of results obtained by the authors on the
geometry of connections with totally skew-symmetric torsion on the
following manifolds: almost complex manifolds with Norden metric,
almost contact manifolds with B-metric and almost hypercomplex
manifolds with Hermitian and anti-Hermitian metric.
\end{abstract}

\maketitle



\section*{Introduction }

In Hermitian geometry there is a strong interest in the
connections preserving the metric and the almost complex structure
whose torsion is totally skew-symmet\-ric
(\cite{Ya,Ya-Kon,Stro,Bi,Gau,Gra-Po,Fri-I,Fri-I2}). Such
connections are called \emph{KT-connections} (or \emph{Bismut
connections}). They find widespread application in mathematics as
well as in theoretic physics. For instance, it is proved a local
index theorem for non-K\"ahler manifolds by KT-connection in
\cite{Bi} and the same connection is applied in string theory in
\cite{Stro}. According to \cite{Gau}, on any Hermitian manifold,
there exists a unique KT-connection. In \cite{Fri-I2} all almost
contact, almost Hermitian and $G_2$-structures admitting a
KT-connection are described.

In this work\footnote{partially supported by projects IS-M-4/2008
and RS09-FMI-003 of the Scientific Research Fund, Paisii
Hilendarski University of Plovdiv, Bulgaria} we provide a survey
of our investigations into connections with totally skew-symmetric
torsion on almost complex manifolds with Norden metric, almost
contact manifolds with B-metric and almost hypercomplex manifolds
with Hermitian and anti-Hermitian metric.



In Section \ref{sec2} we consider an almost complex manifold with
Norden metric (i.e. a neutral metric $g$ with respect to which the
almost complex structure $J$ is an anti-isometry). On such a
manifold we study a natural connection (i.e. a linear connection
$\nabla'$ preserving $J$ and $g$) and having totally
skew-symmetric torsion. We prove that $\nabla'$ exists only when
the manifold belongs to the unique basic class with non-integrable
structure $J$. This is the class $\W_3$ of quasi-K\"ahler
manifolds with Norden metric. We establish conditions for the
corresponding curvature tensor to be K\"ahlerian as well as
conditions $\nabla'$ to have a parallel torsion. We construct a
relevant example on a 4-dimensional Lie group.

In Section \ref{sec3} we consider an almost contact manifold with
B-metric which is the odd-dimensional analogue of an almost
complex manifold with Norden metric. On such a manifold we
introduce the so-called $\ff$KT-connection having totally
skew-symmetric torsion and preserving the almost contact structure
and the metric. We establish the class of the manifolds where this
connection exists. We construct such a connection and study its
geometry. We establish conditions for the corresponding curvature
tensor to be of $\ff$-K\"ahler type as well as conditions for the
connection to have a parallel torsion. We construct an example on
a 5-dimensional Lie group where the $\ff$KT-connection has a
parallel torsion.



In Section \ref{sec4} we consider an almost hypercomplex manifold
with Hermitian and anti-Hermitian metric. This metric is a neutral
metric which is Hermitian with respect to the first almost complex
structure and an anti-Hermitian (i.e. a Norden) metric with
respect to the other two almost complex structures. On such a
manifold we introduce the so-called pHKT-connection having totally
skew-symmetric torsion and preserving the almost hypercomplex
structure and the metric. We establish the class of the manifolds
where this connection exists. We study the unique pHKT-connection
$D$ on a nearly K\"ahler manifold with respect to the first almost
complex structure. We establish that this connection coincides
with the known KT-connection on nearly K\"ahler manifolds and
therefore it has a parallel torsion. We prove the equivalence of
the conditions $D$ be strong, flat and with a parallel torsion
with respect to the Levi-Civita connection.

\section{Almost complex manifold with Norden metric}\label{sec2}

Let $(M,J,g)$ be a $2n$-dimensional \emph{almost complex manifold
with Norden metric}, i.e. $M$ is a differentiable manifold with an
almost complex structure $J$ and a pseudo-Riemannian metric $g$
such that
\[
J^2x=-x, \qquad g(Jx,Jy)=-g(x,y)
\]
for arbitrary $x$, $y$ of the algebra $\X(M)$ on the smooth vector
fields on $M$. Further $x,y,z,w$ will stand for arbitrary elements
of $\X(M)$.

The associated metric $\tilde{g}$ of $g$ on $M$ is defined by
$\tilde{g}(x,y)=g(x,Jy)$. Both metrics are necessarily of
signature $(n,n)$. The manifold $(M,J,\tilde{g})$ is also an
almost complex manifold with Norden metric.

A classification of the almost complex manifolds with Norden
metric is given in \cite{GaBo}. This classification is made with
respect to the tensor $F$
of type (0,3) defined by $
F(x,y,z)=g\bigl( \left( \nabla_x J \right)y,z\bigr)$,
where $\nabla$ is the Levi-Civita connection of $g$. The tensor
$F$ has the following properties
\begin{equation}\label{1.6}
F(x,y,z)=F(x,z,y)=F(x,Jy,Jz),\quad F(x,Jy,z)=-F(x,y,Jz).
\end{equation}
The basic classes are $\W_1$, $\W_2$ and $\W_3$.  Their
intersection is the class $\W_0$ of the K\"ahlerian-type
manifolds, determined by $\W_0: \; F(x,y,z)=0$ $\Leftrightarrow$
$\n J=0.$


The class $\W_3$ of the \emph{quasi-K\"ahler manifolds with Norden
metric} is determin\-ed by the condition
\begin{equation}\label{4}
   \W_3: \quad  F(x,y,z)+F(y,z,x)+F(z,x,y)=0.
\end{equation}

This is the only class of the basic classes $\W_1$, $\W_2$ and
$\W_3$, where each mani\-fold (which is not a K\"ahler-type
manifold) has a non-integrable almost complex structure $J$, i.e.
the Nijenhuis tensor $N$, determined by $N(x,y)=\left(\n_x
J\right)Jy-\left(\n_y J\right)Jx+\left(\n_{Jx}
J\right)y-\left(\n_{Jy}
    J\right)x
$
is non-zero.

The components of the inverse matrix of $g$ are denoted by
$g^{ij}$ with respect to a basis $\{e_i\}$ of the tangent space
$T_pM$ of $M$ at a point $p\in M$.

The \emph{square norm of $\nabla J$} is defined by
$
    \norm{\nabla J}=g^{ij}g^{ks}
    g\bigl(\left(\nabla_{e_i} J\right)e_k,\left(\nabla_{e_j}
    J\right)e_s\bigr).
$

\begin{dfn}[\cite{MekMan-1}]
An almost complex manifold with Norden metric and $\norm{\nabla
J}=0$ is called an \emph{isotropic-K\"ahler manifold}.
\end{dfn}

\subsection{KT-connection}

Let $\nabla'$ be a linear connection on an almost complex manifold
with Norden metric $(M,J,g)$. If $T$ is the torsion tensor of
$\nabla'$, i.e. $T(x,y)=\nabla'_x y-\nabla'_y x-[x, y]$, then the
corresponding tensor of type (0,3) is determined by
$T(x,y,z)=g(T(x,y),z)$.

\begin{dfn}[\cite{Ga-Mi}]
A linear connection $\nabla'$ preserving the almost complex
structure $J$ and the Norden metric $g$, i.e.
$\nabla'J=\nabla'g=0$, is called a \emph{natural connection} on
$(M,J,g)$.
\end{dfn}

By analogy with Hermitian geometry we have given the following
\begin{dfn}[\cite{Mek2}]
A natural connection $\n'$ on an almost complex manifold with
Norden metric is called a \emph{KT-connection} if its torsion
tensor $T$ is totally skew-symmetric, i.e. a 3-form.
\end{dfn}

We have proved the following
\begin{thm}[\cite{Mek7}]
If a KT-connection $\n'$ exists on an almost complex manifold with
Norden metric then the manifold is quasi-K\"ahlerian with Norden
metric.
\end{thm}

A partial decomposition of the space $\mathcal{T}$ of the torsion
(0,3)-tensors $T$ is valid on an almost complex manifold with
Norden metric $(M,J,g)$ according to \cite{Ga-Mi}:
$\mathcal{T}=\mathcal{T}_1\oplus\mathcal{T}_2\oplus\mathcal{T}_3\oplus\mathcal{T}_4$,
where $\mathcal{T}_i$ $(i=1,2,3,4)$ are invariant orthogonal
subspaces.

\begin{thm}[\cite{Mek7}]
Let $\n'$ be a KT-connection with torsion $T$ on a quasi-K\"ahler
manifold with Norden metric $(M,J,g)\notin\W_0$. Then
\begin{itemize}
    \item[1)] $T\in\TT_1\oplus\TT_2\oplus\TT_4$;
    \item[2)] $T$ does not belong to any of the classes $\TT_1\oplus\TT_2$ and $\TT_1\oplus\TT_4$;
    \item[3)] $T\in\TT_2\oplus\TT_4$ if and only if $T$ is determined
    by
\begin{equation}\label{7.2'''}
    T(x,y,z)=-\frac{1}{2}\bigl\{F(x,y,Jz)+F(y,z,Jx)+F(z,x,Jy)\bigr\}.
\end{equation}
\end{itemize}
\end{thm}

Bearing in mind that $T$ is a 3-form,  the following is valid
\begin{equation}\label{7.2}
    g\left(\n'_x y-\n_x y,z\right)=\frac{1}{2}T(x,y,z).
\end{equation}

Then, by \eqref{7.2}, \eqref{1.6} and \eqref{4}, it follows
directly that the tensor $T$, determined by \eqref{7.2'''}, is the
unique torsion tensor of a KT-connection, which is a linear
combination of the components of the basic tensor $F$ on $(M,J,g)$
\cite{Teof}.

Further, the notion of the KT-connection  $\n'$ on $(M,J,g)$ we
refer to the connection with the torsion tensor determined by
\eqref{7.2'''}.

\subsection{KT-connection with K\"ahler curvature tensor or parallel torsion}

\begin{dfn}[\cite{GanGriMih}]\label{dfn-Kaehler}
A tensor $L$ is called a \emph{K\"ahler tensor} if it has the
following properties:
\[
\begin{array}{l}
    L(x,y,z,w)=-L(y,x,z,w)=-L(x,y,w,z),\\[4pt]
    L(x,y,z,w)+L(y,z,x,w)+L(z,x,y,w)=0,\\[4pt]
    L(x,y,Jz,Jw)=-L(x,y,z,w).
\end{array}
\]
\end{dfn}

Let $R'$ be the curvature tensor of the KT-connection $\n'$, i.e.
$R'(x,y)z=\nabla'_x \bigl(\nabla'_y z\bigr) - \nabla'_y
\bigl(\nabla'_x z\bigr) - \nabla'_{[x,y]}z$. The corresponding
tensor of type $(0,4)$ is determined by $R'(x,y,z,w)$
$=g(R'(x,y)z,w)$.

We have therefore proved the following
\begin{thm}[\cite{Mek2}]\label{thm2.1}
The following conditions are equivalent:
\begin{itemize}
    \item[i)]
$R'$ is a K\"ahler tensor;
    \item[ii)]
    $12R'(x,y,z,w)=12R(x,y,z,w)\allowbreak+2g\left(T(x,y),T(z,w)\right)
    \allowbreak-g\left(T(y,z),T(x,w)\right)\allowbreak-g\left(T(z,x),T(y,w)\right)$;
    \item[iii)]   $\mathop{\s} \limits_{x,y,z} \bigl\{
      g\bigl(\left(\nabla_x J\right)Jy+\left(\nabla_{Jx} J\right)y,
      \left(\nabla_z
       J\right)Jw+\left(\nabla_{Jz}J\right)w\bigr)
    \bigr\}=0$,
where $\s$ denotes the cyclic sum by three arguments.
\end{itemize}
\end{thm}

\begin{prop}[\cite{Mek2,Mek6}]
Let $\tau$ and  $\tau'$ be the scalar curvatures for $R$ and $R'$,
respectively. Then the following is valid
\begin{itemize}
    \item[i)]
$3\nJ{}=8(\tau'-\tau)$ if $\n'$ has a K\"ahler curvature tensor;
    \item[ii)]
$\nJ{}=8(\tau-\tau')$ if $\n'$ has a parallel torsion.
\end{itemize}
\end{prop}

\begin{cor}[\cite{Mek6}]
If $\n'$ has a K\"ahler curvature tensor and a parallel torsion
then $(M,J,g)$ is an isotropic-K\"ahler manifold.
\end{cor}



\subsection{An example}

Let $(G,J,g)$ be a 4-dimensional almost complex manifold with
Norden met\-ric, where $G$ is the connected Lie group with an
associated Lie algebra $\mathfrak{g}$ determin\-ed by a global
basis $\{X_i\}$ of left invariant vector fields, and $J$ and $g$
are the almost complex structure and the Norden metric,
respectively, determined by
\[
    JX_1=X_3,\qquad JX_2=X_4,\qquad JX_3=-X_1,\qquad JX_4=-X_2
\]
and
\[
\begin{array}{c}
  g(X_1,X_1)=g(X_2,X_2)=-g(X_3,X_3)=-g(X_4,X_4)=1, \\[4pt]
  g(X_i,X_j)=0\quad \text{for}\quad i\neq j. \\
\end{array}
\]

\begin{thm}[\cite{Mek7}]
The manifold $(G,J,g)$ is a quasi-K\"ahlerian with a Killing
associated Norden metric $\tilde{g}$, i.e.
$g\left([X_i,X_j],JX_k\right)+g\left([X_i,X_k],JX_j\right)=0$, if
and only if $\mathfrak{g}$ is defined by
\[
\begin{array}{ll}
    [X_1,X_2]= \lambda_1 X_1 +\lambda_2 X_2,\qquad & [X_1,X_3]= \lambda_3
    X_2-\lambda_1 X_4,\\[4pt]
    [X_1,X_4]= -\lambda_3 X_1 -\lambda_2 X_4,\qquad & [X_2,X_3]= \lambda_4
    X_2+\lambda_1 X_3,\\[4pt]
    [X_2,X_4]= -\lambda_4 X_1 +\lambda_2 X_3,\qquad & [X_3,X_4]= \lambda_3
    X_3
    +\lambda_4 X_4,\\[4pt]
\end{array}
\]
where $\lambda_1$, $\lambda_2$, $\lambda_3$, $\lambda_4\in\R$.
\end{thm}

Let $(G,J,g)$ be the manifold determined by the conditions in the
last theorem.

The non-trivial components $T_{ijk}=T(X_i,X_j,X_k)$ of the torsion
$T$ of the KT-connection $\n'$ on $(G,J,g)$ are
$T_{134}=\lambda_1$, $T_{234}=\lambda_2$, $T_{123}=-\lambda_3$,
$T_{124}=-\lambda_4$.

Moreover it is proved the following
\begin{thm}[\cite{Mek7}]
The following propositions are equivalent:
\begin{itemize}
    \item[i)] The manifold $(G,J,g)$ is isotropic-K\"ahlerian;
    \item[ii)] The manifold $(G,J,g)$ is scalar flat;
    \item[iii)] The KT-connection $\n'$ has a K\"ahler curvature
    tensor;
    \item[iv)] The equality $\lm_1^2+\lm_2^2-\lm_3^2-\lm_4^2=0$ is valid.
\end{itemize}
\end{thm}

\section{Almost contact manifolds with B-metric}\label{sec3}

Let $(M,\ff,\xi,\eta,g)$ be an \emph{almost contact manifold with
B-metric} (an \emph{almost contact B-metric manifold}), i.e. $M$
is a $(2n+1)$-dimensional differentiable mani\-fold with an almost
contact structure $(\ff,\xi,\eta)$ which consists of an
endomorphism $\ff$ of the tangent bundle, a vector field $\xi$,
its dual 1-form $\eta$ as well as $M$ is equipped with a
pseudo-Riemannian metric $g$  of signature $(n,n+1)$, such that
the following algebraic relations are satisfied
\[
\begin{array}{c}
\ff\xi = 0,\qquad \ff^2 = -I + \eta \otimes \xi,\qquad
\eta\circ\ff=0,\qquad \eta(\xi)=1,\\[4pt]
g(\ff x, \ff y ) = - g(x, y ) + \eta(x)\eta(y),
\end{array}
\]
 where $I$
denotes the identity.

Let us remark that the so-called B-metric $g$ one can say a
\emph{metric of Norden type} in the odd-dimensional case, because
the restriction of $g$ on the contact distribution $\ker \eta$ is
a Norden metric with respect to the almost complex structu\-re
derived by $\ff$.

The associated metric $\tilde{g}$ of $g$ on $M$ is defined by
$\tilde{g}(x,y)=g(x,\ff y)+\eta(x)\eta(y)$. Both metrics are
necessarily of signature $(n,n+1)$. The manifold
$(M,\ff,\xi,\eta,\tilde{g})$ is also an almost contact B-metric
manifold.

A classification of the almost contact manifolds with B-metric is
given in \cite{GaMiGr}. This classification is made with respect
to the tensor $F$ of type (0,3) defined by
$
F(x,y,z)=g\bigl( \left( \nabla_x \ff \right)y,z\bigr)$,
where $\nabla$ is the Levi-Civita connection of $g$. The tensor
$F$ has the following properties
\[
F(x,y,z)=F(x,z,y)=F(x,\ff y,\ff
z)+\eta(y)F(x,\xi,z)+\eta(z)F(x,y,\xi).
\]
This classification includes eleven basic classes $\F_1,\ \F_2,
\dots, \F_{11}$. The special class $\F_0$
, belonging to any other class $\F_i$ $(i=1,2,\dots,11)$, is
determined by the condition $F(x,y,z)=0$. Hence $\F_0$ is the
class of almost contact B-metric manifolds with $\n$-parallel
structures, i.e. $\n\ff=\n\xi=\n\eta=\n g=0$.

In the present work we pay attention to $\F_3$ and $\F_7$, where
each mani\-fold (which is not a $\F_0$-manifold) has a
non-integrable almost contact structure, i.e. the Nijenhuis tensor
$N$, determined by $N(x, y )\allowbreak{} = [\ff, \ff](x, y ) +
d\eta(x, y )\xi$, is non-zero. These basic classes are
characterized by the conditions
\[
\F_3:\quad \sx F(x,y,z)=0,\quad F(\xi,y,z)=F(x,y,\xi)=0,
\]
\[
\F_7:\quad \sx F(x,y,z)=0,\quad F(x,y,z)=-F(\ff x,\ff y,z)-F(\ff
x,y,\ff z).\label{1.4a}
\]

Let us consider the linear projectors $h$ and $v$ over $T_pM$
which split (orthogo\-nal\-ly and invariantly with respect to the
structural group) any vector $x$ into a horizontal component
$h(x)=-\ff^2x$ and a vertical component $v(x)=\eta(x)\xi$.

The decomposition $T_pM=h(T_pM)\oplus v(T_pM)$ generates the
corresponding distribution of basic tensors $F$, which gives the
horizontal component $\F_3$ and the vertical component $\F_7$ of
the class $\F_3\oplus\F_7$.

The \emph{square norm of $\nabla \ff$} is defined by $
    \norm{\nabla \ff}=g^{ij}g^{ks}
    g\bigl(\left(\nabla_{e_i} \ff\right)e_k,\left(\nabla_{e_j}
    \ff\right)e_s\bigr).
$

\begin{dfn}[\cite{Man31}]
An almost contact B-metric manifold with $\norm{\nabla \ff}=0$ is
called an \emph{isotropic-$\F_0$-manifold}.
\end{dfn}

\subsection{$\ff$KT-connection}

\begin{dfn}[\cite{Man31}]
A linear connection $D$ preserving the almost contact B-metric
structure $(\ff,\xi,\eta,g)$, i.e. $D\ff=D\xi=D\eta=Dg=0$, is
called a \emph{natural connection} on $(M,\ff,\xi,\eta,g)$.
\end{dfn}

\begin{dfn}[\cite{Man31}]
A natural connection $D$ on an almost contact B-metric manifold is
called a \emph{$\ff$KT-connection} if its torsion tensor $T$ is
totally skew-symmetric, i.e. a 3-form.
\end{dfn}

The following theorem is proved.
\begin{thm}[\cite{Man31}]
If a $\ff$KT-connection $D$ exists on an almost contact B-metric
manifold $(M,\ff,\xi,\eta,g)$ then $\xi$ is a Killing vector field
and $\s F=0$, i.e. $(M,\ff,\xi,\allowbreak\eta,g)$ belongs to the
class $\F_3\oplus\F_7$.
\end{thm}

The existence of a $\ff$KT-connection $D$ on a manifold in
$\F_3\oplus\F_7$ is given by the following
\begin{prop}[\cite{Man31}]
Let $(M,\ff,\xi,\allowbreak\eta,g)$ be in the class
$\F_3\oplus\F_7$. Then the connection $D$ with a torsion tensor
$T$, determined by
\begin{equation}\label{T37} %
T(x,y,z)=-\frac{1}{2} \sx\bigl\{F(x,y,\ff
z)-3\eta(x)F(y,\ff z,\xi)\bigr\}, %
\end{equation} %
is a $\ff$KT-connection on $(M,\ff,\xi,\eta,g)$.
\end{prop}

Further, the notion of the $\ff$KT-connection  $D$ on
$(M,\ff,\xi,\eta,g)$ we refer to the connection with the torsion
tensor determined by \eqref{T37}. For this connection we have
\[
\begin{split}
D_x y=\n_x y+\frac{1}{4}\bigl\{&2\left(\n_x\ff\right)\ff
y-\left(\n_y\ff\right)\ff x+\left(\n_{\ff y}\ff\right)x \\[4pt]
&+3
\eta(x)\n_y\xi-4\eta(y)\n_x\xi+2\left(\n_x\eta\right)y.\xi\bigr\}.
\end{split}
\]

\subsection{The $\ff$KT-connection on the horizontal component}

Let us consider a manifold from the class $\F_3$ -- the horizontal
component of $\F_3\oplus\F_7$. Since the restriction on the
contact distribution of any $\F_3$-manifold is an almost complex
manifold with Norden metric belonging to the class $\W_3$ (known
as a quasi-K\"ahler manifold with Norden metric), then
    the curvature properties are obtained in a way analogous to that in Section \ref{sec2}.

\subsection{The $\ff$KT-connection on the vertical component}

Let $(M,\ff,\xi,\eta,g)$ belong to the class $\F_7$ -- the
vertical component of $\F_3\oplus\F_7$. For such a manifold the
torsion of
the $\ff$KT-connection $D$ has the form %
\[
T(x,y)=2\left\{\eta(x)\n_y\xi-\eta(y)\n_x\xi+\left(\n_x\eta\right)y.\xi\right\}.
\]

A \emph{tensor of $\ff$-K\"ahler type} we call a tensor with the
properties from \dfnref{dfn-Kaehler} with respect to the structure
$\ff$.

We have proved the following
\begin{thm}[\cite{Man31}]
The curvature tensor $K$ of  $D$ on a $\F_7$-ma\-ni\-fold is of
$\ff$-K\"ahler type if and only if it has the form
\[
\begin{aligned}
K(x,y,z,w) = R(x, y, z, w) %
+\frac{1}{3}\bigl\{%
2 \left(\n_x\eta\right)y\left(\n_z\eta\right)w%
-\left(\n_y\eta\right)z\left(\n_x \eta\right)w\\[4pt]%
-\left(\n_z \eta\right)x\left(\n_y \eta\right)w\bigr\}\\[4pt]
+\eta(x)\eta(z)g\left(\n_y \xi,\n_w \xi\right)%
-\eta(x)\eta(w)g\left(\n_y \xi,\n_z \xi\right)\phantom{.}\\[4pt]
-\eta(y)\eta(z)g\left(\n_x \xi,\n_w \xi\right)%
+\eta(y)\eta(w)g\left(\n_x \xi,\n_z \xi\right).
\end{aligned}
\]
\end{thm}

\begin{thm}[\cite{Man31}]
If $D$ has a curvature tensor $K$ of $\ff$-K\"ahler type and a
parallel torsion $T$ on a $\F_7$-manifold then
\[
\begin{aligned}
K(x,y,z,w) = R(x, y, z, w) +\frac{1}{3}\bigl\{2 \left(\n_x
\eta\right)y\left(\n_z
\eta\right)w%
&+\left(\n_x \eta\right)z\left(\n_y \eta\right)w\\[4pt]
&-\left(\n_x \eta\right)w\left(\n_y \eta\right)z\bigr\},
\end{aligned}
\] %
\[
\rho(K)(y,z)=\rho(y,z), \qquad \tau(K)=\tau,
\]
where $\rho(K)$ and $\rho$ are the Ricci tensors for $K$ and $R$,
respectively, and $\tau(K)$ and $\tau$ are the their corresponding
scalar curvatures.

\end{thm}

\subsection{An example}

Let $(G,\ff,\xi,\eta,g)$ be a 5-dimensional almost contact
manifold with B-met\-ric, where $G$ is the connected Lie group
with an associated Lie algebra $\mathfrak{g}$ determin\-ed by a
global basis $\{X_i\}$ of left invariant vector fields, and
$(\ff,\xi,\eta)$ and $g$ are the almost contact structure and the
B-metric, respectively, determined by
\[
\begin{array}{c}
\ff X_1 = X_3,\quad \ff X_2 = X_4,\quad \ff X_3 =-X_1,\quad \ff
X_4 =
-X_2,\quad \ff X_5 =0;\\[4pt]
\xi=X_5;\qquad \eta(X_i)=0\; (i=1,2,3,4),\quad \eta(X_5)=1;
\end{array}
\]
\[
\begin{array}{c}
g(X_1,X_1)=g(X_2,X_2)=-g(X_3,X_3)=-g(X_4,X_4)=g(X_5,X_5)=1,
\\[4pt]
g(X_i,X_j)=0,\; i\neq j,\quad  i,j\in\{1,2,3,4,5\}.
\end{array}
\]

\begin{thm}[\cite{Man31}]
The manifold $(G,\ff,\xi,\eta,g)$ is a $\F_7$-manifold if and only
if $\mathfrak{g}$ is determined by the following non-zero
commutators:
\[
\begin{array}{c}
[X_1,X_2]=-[X_3,X_4]=-\lm_1X_1-\lm_2X_2+\lm_3X_3+\lm_4X_4+2\mu_1X_5,
\\[4pt]
[X_1,X_4]=-[X_2,X_3]=-\lm_3X_1-\lm_4X_2-\lm_1X_3-\lm_2X_4+2\mu_2X_5,
\end{array}
\]
where $\lm_i, \mu_j \in \R$ $(i=1,2,3,4; j=1,2)$.
\end{thm}

Let $(G,\ff,\xi,\eta,g)$ be the manifold determined by the
conditions in the last theorem.

The non-trivial components $T_{ijk}=T(X_i,X_j,X_k)$ of the torsion
$T$ of the $\ff$KT-connection $D$ on $(G,\ff,\xi,\eta,g)$ are
$T_{125}=T_{345}=2\mu_1$, $T_{235}=T_{415}=2\mu_2$.

Hence, using the components of $D$, we calculate that the
corresponding compo\-nents of the covariant derivative of $T$ with
respect to $D$ are zero. Thus, we have proved the following
\begin{thm}[\cite{Man31}]
The $\ff$KT-connection $D$ on $\Lf$ has a parallel torsion $T$.
\end{thm}

\begin{thm}[\cite{Man31}]
The manifold $\Lf$ is an isotropic-$\F_0$-manifold if and only if
$\mu_1=\pm\mu_2$.
\end{thm}

\section{Almost hypercomplex manifolds with Hermitian and Norden met\-ric}\label{sec4}

Let $(M,H)$ be an \emph{almost hypercomplex manifold}, i.e. $M$ is
a $4n$-dimension\-al differentiable manifold and $H=(J_1,J_2,J_3)$
is a triple of almost complex structures with the properties:
\[
J_\al=J_\bt\circ J_\gm=-J_\gm\circ J_\bt, \qquad J_\al^2=-I%
\]
for all cyclic permutations $(\al, \bt, \gm)$ of $(1,2,3)$.

The standard structure of $H$ on a $4n$-dimensional vector space
with a basis\\
$\{X_{4k+1},X_{4k+2},X_{4k+3},X_{4k+4}\}_{k=0,1,\dots,n-1}$ has
the form:
\[
\begin{array}{lll}
J_1X_{4k+1}= X_{4k+2},\;\; & J_2X_{4k+1}=X_{4k+3},\;\; & J_3X_{4k+1}=-X_{4k+4},\\[4pt]
J_1X_{4k+2}=-X_{4k+1},\;\; & J_2X_{4k+2}=X_{4k+4}, \;\; & J_3X_{4k+2}=X_{4k+3}, \\[4pt]
J_1X_{4k+3}=-X_{4k+4},\;\; & J_2X_{4k+3}=-X_{4k+1}, \;\; & J_3X_{4k+3}=-X_{4k+2},\\[4pt]
J_1X_{4k+4}=X_{4k+3},\;\; & J_2X_{4k+4}=-X_{4k+2}, \;\; &
J_3X_{4k+4}=X_{4k+1}.
\end{array}
\]

Let $g$ be a pseudo-Riemannian metric on $(M,H)$ with the
properties
\[
g(x,y)=\ea g(J_\al x,J_\al y), \qquad \ea=
\begin{cases}
\begin{array}{ll}
1, \quad & \al=1;\\
-1, \quad & \al=2;3.
\end{array}
\end{cases}
\]
In other words, for $\al=1$, the metric $g$ is Hermitian with
respect to $J_1$, whereas in the cases  $\al=2$ and $\al=3$ the
metric $g$ is an anti-Hermitian (i.e. Norden) metric with respect
to $J_2$ and $J_3$, respectively. Moreover, the associated
bilinear forms $g_1$, $g_2$, $g_3$ are determined by
\[
g_\al(x,y)=g(J_\al x,y)=-\ea g(x,J_\al y),\qquad
\al=1,2,3. %
\]

Then, we call a manifold with such a structure briefly an
\emph{almost $(H,G)$-manifold} \cite{GrMa,GrMaDi}.

The structural tensors of an almost $(H,G)$-manifold are the three
$(0,3)$-tensors determined by
\[
F_\al (x,y,z)=g\bigl( \left( \n_x J_\al
\right)y,z\bigr)=\bigl(\n_x g_\al\bigr) \left( y,z \right),\qquad
\al=1,2,3,
\]
where $\n$ is the Levi-Civita connection generated by $g$.

In the classification of Gray-Hervella \cite{GrHe} for almost
Hermitian manifolds the class $\G_1=\W_1\oplus\W_3\oplus\W_4$ is
determined by the condition $F_1(x,x,z)=F_1(J_1x,J_1x,z)$.

\begin{thm}[\cite{Man32}]
If $M$ is an almost $(H,G)$-manifold which is a quasi-K\"ahler
manifold with Norden metric regarding $J_2$ and $J_3$, then it
belongs to the class $\G_1$ with respect to $J_1$.
\end{thm}

\subsection{pHKT-connection}

\begin{dfn}[\cite{Man32}]
A linear connection $D$ preserving the almost hypercomplex
structure $H$ and the metric $g$, i.e. $DJ_1=DJ_2=DJ_3=Dg=0$, is
called a \emph{natural connection} on $(M,H,G)$.
\end{dfn}

\begin{dfn}[\cite{Man32}]
A natural connection $D$ on an almost $(H,G)$-manifold is called a
\emph{pseudo-HKT-connection} (briefly, a \emph{pHKT-connec\-tion})
if its torsion tensor $T$ is totally skew-symmetric, i.e. a
3-form.
\end{dfn}

For an almost complex manifold with Hermitian metric $(M,J,g)$, in
\cite{Fri-I2} it is proved that there exists a unique
KT-connection if and only if the Nijenhuis tensor $N_{J}(x,y,z):=
g(N_{J}(x,y),z)$ is a 3-form, i.e. the manifold belongs to the
class of cocalibrated structures $\G_1$.

\subsection{The class $\W_{133}$}

Next, we restrict the class $\G_1(J_1)$ to its subclass
$\W_1(J_1)$  of \emph{nearly K\"ahler manifolds with neutral
metric} regarding $J_1$ defined by $F_1(x,y,z)=-F_1(y,x,z)$. In
this case $(M,H,G)$ belongs to the class
$\W_{133}=\W_1(J_1)\cap\W_3(J_2)\cap\W_3(J_3)$ and $\dim{M}\geq
8$.

We have proved the following
\begin{thm}[\cite{Man32}]
The curvature tensor $R$ of $\n$ on $(M,H,G)\in\W_{133}$ has the
following property with respect to the almost hypercomplex
structure $H$:
\[
    R(x,y,z,w)+\sum_{\al=1}^3 R(x,y,J_\al z,J_\al w)
    =\sum_{\al=1}^3\bigl\{A_\al(x,z,y,w)-A_\al(y,z,x,w)\bigr\},
\]
where \(
A_\al(x,y,z,w)=g\bigl(\left(\n_xJ_\al\right)y,\left(\n_zJ_\al\right)w\bigr)
\), $\al=1,2,3$.
\end{thm}

\subsection{The pHKT-connection on a $\W_{133}$-manifold}

KT-connections on nearly K\"ahler manifolds are investigated for
instance in \cite{Nagy}. The unique KT-connection $D^1$ for the
nearly K\"ahler manifold $(M,J_1,g)$ on the considered almost
$(H,G)$-manifold has the form
\[
    g\left(D^1_xy,z\right)=g\left(\n_xy,z\right)+\frac{1}{2}F_1(x,y,J_1z).
\]

Moreover, there exists a unique KT-connection $D^\al$ ($\al=2,3$)
for the quasi-K\"ahler manifold with Norden metric $(M,J_\al,g)$
on the considered almost $(H,G)$-manifold such that
\[
    g\left(D^\al_xy,z\right)=g\left(\n_xy,z\right)-\frac{1}{4}\sx F_\al(x,y,J_\al z).
\]

In \cite{Man32} we have constructed a connection $D$, using the
KT-connections $D^1$, $D^2$ and $D^3$, on an almost
$(H,G)$-manifold from the class $\W_{133}$ and we have proved the
following
\begin{thm}[\cite{Man32}]\label{prop-HKT}
The connection $D$ defined by
\[
    g\left(D_xy,z\right)=g\left(\n_xy,z\right)+\frac{1}{2}F_1(x,y,J_1z).
\]
is the unique pHKT-connection on an almost $(H,G)$-manifold from
the class $\W_{133}$.
\end{thm}

Let us remark that the pHKT-connection $D$ on an almost
$(H,G)$-manifold coincides with the known KT-connection $D^1$ on
the corresponding nearly K\"ahler manifold. Then the torsion of
the pHKT-connection $D$ is parallel and henceforth $T$ is
coclosed, i.e. $\delta T=0$ \cite{Fri-I2}. Moreover, the
cur\-va\-ture tensors $K$ of $D$ and $R$ of $\n$ has the following
relation \cite{Gray}
\[
    K(x,y,z,w)=R(x,y,z,w)+\frac{1}{4} A_1(x,y,z,w)+\frac{1}{4}\sx A_1(x,y,z,w).
\]

We have proved the following
\begin{thm}[\cite{Man32}]\label{prop-sA}
Let $(M,H,G)$ be an almost $(H,G)$-manifold from $\W_{133}$ and
$D$ be the pHKT-connection. Then the following characte\-ristics
of this connection are equivalent:
\begin{enumerate}\renewcommand{\labelenumi}{(\roman{enumi})}
    \item $D$ is strong ($\dd T=0$);
    \item $D$ has a $\n$-parallel torsion;
    \item $D$ is flat.
\end{enumerate}
\end{thm}

\begin{thm}[\cite{Man32}]\label{prop-flat}
Let $(M,H,G)$ be an almost $(H,G)$-manifold from $\W_{133}$ and
$D$ be the pHKT-connection. If $D$ is flat or strong then
$(M,H,G)$ is $\n$-flat, isotropic-hyper-K\"ahlerian (i.e.
$\nJ{\al}=0$, $\al=1,2,3$) and the torsion of $D$ is isotropic
(i.e. $\nT=0$).
\end{thm}

\end{document}